\documentclass[11pt,oneside,reqno]{amsart}
\usepackage{bm,cancel}
\IfFileExists{ulem.sty}{\usepackage{ulem}}{}
\usepackage{amssymb,amsfonts}
\usepackage{amsmath,latexsym,enumerate}
\usepackage{mathrsfs}

\usepackage[top=1in, bottom=1in, left=1.25in, right=1.25in, marginparwidth=1.1in, marginparsep=0.05in]{geometry}

\usepackage[unicode,breaklinks=true,colorlinks=true]{hyperref}

\usepackage[svgnames]{xcolor}

\newtheorem{thm}{Theorem}[section]
\newtheorem{cor}[thm]{Corollary}
\newtheorem{lem}[thm]{Lemma}
\newtheorem{prop}[thm]{Proposition}
\theoremstyle{definition}

\newtheorem{rem}[thm]{Remark}
\numberwithin{equation}{section}
\newcommand{\R}{\mathbb R}

\newcommand{\cF}{\mathcal{F}}

\usepackage{accents}
\newcommand*{\dotc}[1]{\accentset{\mbox{\large\bfseries .}}{#1}}

\renewcommand{\div}{\mathop{\rm div}\nolimits}

\newcommand{\C}{\mathbb C}
\newcommand{\E}{\mathbb E}
\newcommand{\cL}{\mathcal L}
\newcommand{\cB}{\mathcal B}
\newcommand{\cX}{\mathcal X}
\newcommand{\cY}{\mathcal Y}
\newcommand{\cR}{\mathcal R}

\newcommand{\half}{\R^d_+}
\newcommand{\lN}{\ell^2_N}
\newcommand{\pair}[2]{\langle #1,#2\rangle}
\DeclareMathOperator{\Tr}{Tr}
\DeclareMathOperator{\supp}{supp}
\DeclareMathOperator{\Rea}{Re}
\allowdisplaybreaks[2]
\begin{document}
\title{Maximal regularity for the Stokes operator in bounded $C^1$ domains}
\author[S. Liang]{Su Liang}
\address[S. Liang]{Institute for Theoretical Sciences, Westlake University,
Hangzhou, Zhejiang 310030, China}
\email{liangsu@westlake.edu.cn}
\date{}
\begin{abstract}
Let $\Omega$ be a bounded $C^1$ domain in $\R^d$, $d\geq2$.
We establish vector-valued resolvent estimates for the
Dirichlet Stokes system in $L^q(\Omega;\ell^2)$ for every
$1<q<\infty$.
The estimates are uniform in the number and choice of
spectral parameters in any fixed sector avoiding the
negative real axis.
These estimates imply that the Stokes operator is
$\mathcal R$-sectorial of angle zero and has maximal
$L^r$-regularity for every $1<r<\infty$.
Our proof builds on the perturbation approach of Geng and
Shen, with a different half-space argument to handle
independently varying spectral parameters.
\end{abstract}
\maketitle

\tableofcontents

\section{Introduction}

Let $\Omega$ be a bounded connected $C^1$ domain in $\R^d$,
$d\geq2$, and put
\[
 \Sigma_\theta=\{z\in\C\setminus\{0\}:|\arg z|<\pi-\theta\},
 \qquad 0<\theta<\pi/2.
\]
We consider families of Dirichlet Stokes problems
\begin{equation}\label{eq:S-resolvent}
 \left\{\begin{aligned}
 -\Delta u_j+\nabla p_j+\lambda_j u_j&=F_j&&\text{in }\Omega,\\
 \div u_j&=0&&\text{in }\Omega,\\
 u_j&=0&&\text{on }\partial\Omega,
 \end{aligned}\right.
 \qquad \lambda_j\in\Sigma_\theta,
\end{equation}
where $j=1,\ldots,N$ and $N\in\mathbb N\cup\{\infty\}$.
We establish estimates in $L^q(\Omega;\lN)$ that are
uniform in the number and choice of the parameters.
These estimates yield maximal regularity for the
nonstationary Stokes problem with independent spatial
and temporal exponents.

Resolvent estimates and maximal regularity for the
Stokes operator have been studied extensively on smooth
domains. Classical contributions include the maximal
regularity estimates of Solonnikov \cite{Solonnikov-1977},
the analyticity result of Giga \cite{Giga-1981}, and
the estimates with independent spatial and temporal
exponents of Giga and Sohr \cite{Giga-Sohr-1991}.
Farwig and Sohr \cite{Farwig-Sohr-1994}
developed generalized resolvent estimates and a
perturbation approach for bounded and unbounded domains.
Desch, Hieber and Pr\"uss \cite{Desch-Hieber-Pruss-2001}
established maximal regularity for the half-space problem.
Weis \cite{Weis-2001} characterized maximal regularity
in terms of $\cR$-boundedness of the resolvent through
an operator-valued Fourier multiplier theorem;
see also \cite{Denk-Hieber-Pruss-2003}.
Geissert, Hess, Hieber, Schwarz and Stavrakidis
\cite{Geissert-et-al-2010} proved maximal regularity
and $\cR$-sectoriality by localization using
Bogovskii's operator.
Abe, Giga and Hieber \cite{Abe-Giga-Hieber-2015}
obtained $L^\infty$ resolvent estimates on smooth domains.
For an overview of the smooth-domain theory, see
\cite{Hieber-Saal-2018}.

On domains with nonsmooth boundaries, the admissible
exponents depend on the boundary regularity.
For stationary boundary value problems and Sobolev
regularity on Lipschitz and $C^1$ domains, see
\cite{Fabes-Kenig-Verchota-1988,Brown-Shen-1995,Dindos-Mitrea-2004,Mitrea-Monniaux-2008,Mitrea-Wright-2012}.
For bounded Lipschitz domains in dimension $d\geq3$,
Shen \cite{Shen-2012} established resolvent estimates when
\[
 \left|\frac1q-\frac12\right|<\frac1{2d}+\varepsilon,
\]
where $\varepsilon>0$ depends on the domain.
Extending Shen's approach, Kunstmann and Weis
\cite{Kunstmann-Weis-2017} established maximal regularity
estimates in the same range of $q$.
In this range, Tolksdorf \cite{Tolksdorf-2018}
identified the domain of the square root of the Stokes operator
and obtained gradient and semigroup estimates.
Gabel and Tolksdorf \cite{Gabel-Tolksdorf-2022}
obtained resolvent and maximal regularity estimates
on bounded Lipschitz domains in dimension two for
$|1/q-1/2|<1/4+\varepsilon$.
For bounded $C^1$ domains in any dimension $d\geq2$,
Geng and Shen \cite{Geng-Shen-2025} proved the
single-parameter estimate corresponding to
\eqref{eq:main-est} below for every $1<q<\infty$,
and consequently obtained the analyticity of the
Stokes semigroup.
See also \cite{Monniaux-Shen-2018,Kunstmann-Tolksdorf-2025} for related
results and further references.

Under Sobolev-multiplier conditions on the boundary
charts, Breit \cite{Breit-2025} established maximal
regularity estimates with separate control of
$\partial_tu$, $\nabla^2u$ and $\nabla p$.
Breit and Gaudin \cite{Breit-Gaudin-2025} developed
a Sobolev and Besov regularity theory for the Stokes
equations, including maximal regularity estimates,
under related multiplier assumptions.
The latter results apply, in particular, to bounded
$C^{1,\alpha}$ domains for every $0<\alpha\leq1$,
but the multiplier conditions are not automatic
for a general $C^1$ boundary.
At $q=\infty$, Geng and Shen \cite{Geng-Shen-2026}
established resolvent estimates for the velocity
on bounded $C^1$ domains in dimension $d\geq3$
and bounded Lipschitz domains in dimension $d=2$.
More recently, Shen \cite{Shen-2026} proved
$L^q$ resolvent estimates on three-dimensional
bounded Lipschitz domains for
$(3/2)-\varepsilon<q\leq\infty$.
These last two works assume that the force $F$
is solenoidal, namely $\div F=0$ in $\Omega$
and $F\cdot n=0$ on $\partial\Omega$ in the weak sense.
For $1<q<\infty$, resolvent estimates for arbitrary
$L^q$ forces require boundedness of the Helmholtz
projection on $L^q$, which may fail for both small
and large $q$ on general Lipschitz domains;
see \cite[Appendix B]{Shen-2026}.
At $q=1$, Binz and Gaudin \cite{Binz-Gaudin-2026}
obtained an analytic Stokes semigroup on the quotient
of $L^1$ by gradient fields for bounded $C^{1,\alpha}$
domains, whereas the solenoidal realization does not
generate a $C_0$-semigroup.
Here we establish vector-valued resolvent estimates
and maximal regularity on bounded $C^1$ domains,
without additional multiplier conditions on the boundary.

\begin{samepage}
Our principal result is the following vector-valued
resolvent estimate.\nopagebreak[4]

\begin{thm}\label{thm:main}
Let $\Omega$ be a bounded connected $C^1$ domain in $\R^d$,
$d\geq2$, and let $1<q<\infty$.
For $\lambda_j\in\Sigma_\theta$, put
$\rho_j=(1+|\lambda_j|)^{1/2}$.
For every $(F_j)\in L^q(\Omega;\lN(\C^d))$,
the problems \eqref{eq:S-resolvent} have unique solutions
$u_j\in W^{1,q}_0(\Omega;\C^d)$ and
$p_j\in L^q_0(\Omega)$ satisfying
\begin{equation}\label{eq:main-est}
 \|(\rho_j\nabla u_j)_j\|_{L^q(\Omega;\lN)}
 +\|(\rho_j^2u_j)_j\|_{L^q(\Omega;\lN)}
 \leq C\|(F_j)_j\|_{L^q(\Omega;\lN)}.
\end{equation}
The constant $C$ depends only on $d,q,\theta$ and $\Omega$.
\end{thm}
\end{samepage}

The parameters $\lambda_j$ may vary independently.
For $q\neq2$, uniform bounds for the individual
resolvents do not by themselves imply
\eqref{eq:main-est}.
The vector-valued estimate gives $\cR$-boundedness
of both the resolvent family and its weighted gradient
family.
Section \ref{sec:bounded} also gives estimates for
the pressure, divergence-form forces and nonzero
prescribed divergence.

We next state the maximal regularity consequences.
Let
\[
 L^q_\sigma(\Omega)=
 \overline{\{v\in C^\infty_c(\Omega;\C^d):
                    \div v=0\}}^{\,L^q}.
\]
On a bounded $C^1$ domain, the Helmholtz decomposition
holds for every $1<q<\infty$;
see \cite{Fabes-Mendez-Mitrea-1998} and
\cite[Chapter III]{Galdi-2011}.
We use the Stokes realization of
\cite{Geng-Shen-2025}, defined on $L^q_\sigma(\Omega)$ by
\begin{equation}\label{eq:Stokes-domain}
 \begin{split}
 D(A_q)&=\{u\in W^{1,q}_0(\Omega;\C^d)\cap L^q_\sigma(\Omega):
 -\Delta u+\nabla p\in L^q_\sigma(\Omega) \text{ for some }p\in L^q_0(\Omega)\},\\[-2pt]
 A_qu&=-\Delta u+\nabla p.
 \end{split}
\end{equation}
The derivatives are understood in distributions.

\begin{cor}\label{thm:MR}
Under the assumptions of Theorem \ref{thm:main},
$A_q$ is invertible and $\cR$-sectorial of angle zero
on $L^q_\sigma(\Omega)$.
For $1<r<\infty$, $0<T\leq\infty$ and
$f\in L^r(0,T;L^q_\sigma(\Omega))$, there is a unique
\[
 u\in W^{1,r}(0,T;L^q_\sigma(\Omega))
       \cap L^r(0,T;D(A_q))
\]
satisfying $u'+A_qu=f$ and $u(0)=0$.
Moreover,
\begin{equation}\label{eq:MR}
 \|u'\|_{L^r(0,T;L^q(\Omega))}
 +\|A_qu\|_{L^r(0,T;L^q(\Omega))}
 +\|u\|_{L^r(0,T;W^{1,q}(\Omega))}
 \leq C\|f\|_{L^r(0,T;L^q(\Omega))},
\end{equation}
where $C$ is independent of $T$.
\end{cor}

By Theorem \ref{thm:main} and the square-function
characterization of $\cR$-boundedness
\cite{Denk-Hieber-Pruss-2003}, the family
\[
 \{\lambda(\lambda+A_q)^{-1}:\lambda\in\Sigma_\theta\}
\]
is $\cR$-bounded for every $0<\theta<\pi/2$.
Thus $A_q$ is $\cR$-sectorial of angle zero.
Weis's theorem \cite[Theorem 4.2]{Weis-2001}
on the UMD space $L^q_\sigma(\Omega)$ gives the
$L^r(0,T;L^q(\Omega))$ estimates for $u'$ and $A_qu$
in \eqref{eq:MR}, with a constant independent of $T$.
The zero-parameter estimate in Theorem \ref{thm:general}
also gives invertibility of $A_q$ and
\[
 \|v\|_{W^{1,q}(\Omega)}
 \lesssim \|A_qv\|_{L^q(\Omega)},
 \qquad v\in D(A_q).
\]
Applying this estimate to $u(t)$ and integrating in time
gives the $W^{1,q}$ bound in \eqref{eq:MR}.
This proves Corollary \ref{thm:MR}.

For general $F\in L^r(0,T;L^q(\Omega;\C^d))$,
apply Corollary \ref{thm:MR} to $P_qF$, where $P_q$
denotes the Helmholtz projection.
The normalized pressure is recovered from the stationary
problem with datum $F-\partial_tu$ using
Theorem \ref{thm:general} at $\lambda=0$, which also
yields its $L^r(0,T;L^q(\Omega))$ estimate.

\begin{rem}
Estimate \eqref{eq:MR} controls the combined expression
$-\Delta u+\nabla p$ in $L^r(0,T;L^q(\Omega))$.
Global $W^{2,q}$ regularity need not hold on
a general $C^1$ domain.
\end{rem}

Our proof follows the whole-space, half-space,
graph perturbation and localization scheme of
Geng and Shen \cite{Geng-Shen-2025}.
In the whole space, their argument extends directly
using the vector-valued Mikhlin theorem.
In the half-space, we use a different argument to
handle the independent parameters $\lambda_j$
in a single $L^q(\ell^2)$ estimate.
Their scaling reduction to $|\lambda|=1$ is not
available for a general parameter family.
We instead apply the vector-valued Mikhlin theorem
in the tangential variables at each fixed height,
and obtain bounds that also control the subsequent
integration in the normal variable in $L^q(\ell^2)$.

An important feature of their approach is the
compatible choice of the solution and data spaces
$X_\psi^q$ and $Y_\psi^q$.
Their elegant proof that $S_\psi^\lambda$ is a bounded
isomorphism between these spaces \cite[Theorem 4.3]{Geng-Shen-2025} provides the framework
for the graph perturbation.
Our choice of the parameter-dependent solution and data
spaces $\cX_{\Lambda,\Omega}$ and $\cY_{\Lambda,\Omega}$,
defined in \eqref{eq:Xspace} and \eqref{eq:Yspace},
is inspired by this formulation.

The paper is organized as follows.
Section \ref{sec:prelim} records the function spaces
and vector-valued tools.
Section \ref{sec:model} treats the whole space
and half-space.
Section \ref{sec:graph} proves the small Lipschitz
graph estimate.
Section \ref{sec:bounded} establishes the generalized
bounded-domain estimate and proves Theorem \ref{thm:main}.

\section{Notation and preliminaries}\label{sec:prelim}

Throughout the paper, $N\in\mathbb N\cup\{\infty\}$, and the index
$j$ ranges over $\{1,\ldots,N\}$.
Estimates for countable sequences follow by truncation from
the corresponding finite-sequence estimates with constants
independent of $N$.
For a sequence $z=(z_j)_{j=1}^N$ of scalar, vector, or matrix fields,
\[
 \|z\|_{L^q(\Omega;\lN)}
 =\left\|\left(\sum_{j=1}^N|z_j|^2\right)^{1/2}\right\|_{L^q(\Omega)}.
\]
The absolute value denotes the Euclidean or Frobenius norm, as
appropriate. 
For convenience, we sometimes write
$\|X,Y\|=\|X\|+\|Y\|$, with the same convention for more terms.

We write $X\lesssim_\sigma Y$ if $X\leq C(\sigma)Y$,
and omit the subscript when the dependence is clear.
Unless otherwise indicated, implicit constants may depend on
$d,q,\theta,\Omega$,
but not on $N$, $\Lambda$, or the functions being estimated.
The notation $X\simeq_\sigma Y$ means both
$X\lesssim_\sigma Y$ and $Y\lesssim_\sigma X$.

For $\Lambda=(\lambda_j)_{j=1}^N$, put $r_j=|\lambda_j|^{1/2}$.
Multiplication by a scalar sequence is understood componentwise;
for instance, we write $r^2u=(r_j^2u_j)_j$ and use both notations interchangeably.
For spatial coordinates $x=(x_1,\ldots,x_d)$ and vector
fields $H=(H_1,\ldots,H_d)$, write
$x'=(x_1,\ldots,x_{d-1})$,
$H'=(H_1,\ldots,H_{d-1})$, and
$\nabla'=(\partial_1,\ldots,\partial_{d-1})$.
Our conventions for the gradient and the divergence of a matrix
field are
\begin{equation*}
 (\nabla u)_{\alpha k}=\partial_k u_\alpha,\qquad
 (\div f)_\alpha=\sum_{k=1}^d\partial_k f_{\alpha k}.
\end{equation*}
In particular, $\div(\nabla u)=\Delta u$. Repeated spatial
indices are summed; the sequence index $j$ is summed only
when indicated.

Our Fourier convention in dimension $m$ is
\[
 \widehat h(\xi)=(2\pi)^{-m/2}
      \int_{\R^m}e^{-ix\cdot\xi}h(x)\,dx .
\]
We use the same notation for the full Fourier transform
($m=d$) and the partial Fourier transform in
$x'\in\R^{d-1}$ ($m=d-1$), when the variables make the
meaning clear.

\subsection{Function spaces}

Let $1<q<\infty$ and $q'=q/(q-1)$.
When $|\Omega|<\infty$, set
\[
 L^q_0(\Omega)=
 \left\{h\in L^q(\Omega):\int_\Omega h=0\right\}.
\]
The space $W^{1,q}(\Omega;\lN)$ consists of sequences
$u=(u_j)_j$, with $u_j\in W^{1,q}(\Omega)$, such that
\begin{equation*}
 \|u\|_{W^{1,q}(\Omega;\lN)}
 =\|(u_j)_j\|_{L^q(\Omega;\lN)}
       +\|(\nabla u_j)_j\|_{L^q(\Omega;\lN)}<\infty.
\end{equation*}
The subspace $W^{1,q}_0(\Omega;\lN)$ consists of those
$u\in W^{1,q}(\Omega;\lN)$ for which
$u_j\in W^{1,q}_0(\Omega)$ for every $j$.

We denote by $\dotc W^{1,q}(\Omega)$ the homogeneous Sobolev
space modulo constants, equipped with the norm
\[
 \|u\|_{\dotc W^{1,q}(\Omega)}
       =\|\nabla u\|_{L^q(\Omega)}.
\]
Its negative-order counterpart is defined by
\[
 \dotc W^{-1,q}(\Omega)
       =\big(\dotc W^{1,q'}(\Omega)\big)^*.
\]
We write $\pair{g}{H}$ for the duality pairing.
By the Hahn--Banach theorem, each
$g\in\dotc W^{-1,q}(\Omega)$ admits a divergence
representation $g=\div G$, in the sense that
\[
 \pair{g}{H}=-\int_\Omega G\cdot\nabla H,\qquad
 \|g\|_{\dotc W^{-1,q}(\Omega)}
       =\inf_{\div G=g}\|G\|_{L^q(\Omega)}.
\]
For sequences modulo constant sequences, we define
\[
 \|u\|_{\dotc W^{1,q}(\Omega;\lN)}
 =\|(\nabla u_j)_j\|_{L^q(\Omega;\lN)}.
\]
The corresponding dual space is
\begin{equation*}
 \begin{split}
 \dotc W^{-1,q}(\Omega;\lN)
   &=\big(\dotc W^{1,q'}(\Omega;\lN)\big)^*,\\
 \|g\|_{\dotc W^{-1,q}(\Omega;\lN)}
   &=\sup_{\|H\|_{\dotc W^{1,q'}(\Omega;\lN)}\leq1}
       \left|\sum_{j=1}^N\pair{g_j}{H_j}\right|.
 \end{split}
\end{equation*}

Following \cite[Remark 5.3]{Geng-Shen-2025} and the above convention,
we define
\begin{equation}\label{eq:two-negative}
 \begin{split}
 W^{-1,q}(\Omega;\lN)
   &=\big(W^{1,q'}_0(\Omega;\lN)\big)^*,\\
 W^{-1,q}_0(\Omega;\lN)
   &=\big(W^{1,q'}(\Omega;\lN)\big)^*.
 \end{split}
\end{equation}
For $g\in\dotc W^{-1,q}(\Omega;\lN)$, the definitions
immediately give
\begin{equation*}
 \|g\|_{W^{-1,q}(\Omega;\lN)}
 \leq\|g\|_{W^{-1,q}_0(\Omega;\lN)}
 \leq\|g\|_{\dotc W^{-1,q}(\Omega;\lN)}.
\end{equation*}

\subsection{Randomization and vector-valued tools}

The following Khintchine inequality is a standard tool for establishing vector-valued inequalities from their scalar counterparts. See \cite[Proposition 4.5]{Wolff-2003} for a proof.

\begin{lem}[Khintchine inequality]\label{lem:Khintchine}
Let $\varepsilon_j$ be independent random variables taking
the values $1$ and $-1$ with equal probability.
For $0<q<\infty$ and $(a_j)_j\in\lN$,
\begin{equation*}
 \left(\E\left|\sum_{j=1}^N\varepsilon_ja_j\right|^q
                                                    \right)^{1/q}
 \simeq_q\left(\sum_{j=1}^N|a_j|^2\right)^{1/2}.
\end{equation*}
The constants depend only on $q$.
\end{lem}

We next collect some standard Poincaré and Sobolev inequalities in the vector-valued setting.

\begin{lem}[Poincar\'e and Sobolev estimates]\label{lem:Sobolev}
Let $\Omega$ be a bounded connected Lipschitz domain.
For $1<s<\infty$ and $H\in\dotc W^{1,s}(\Omega;\lN)$,
\begin{equation}\label{eq:Poincare}
 \|H-H_\Omega\|_{W^{1,s}(\Omega;\lN)}
 \lesssim \|\nabla H\|_{L^s(\Omega;\lN)},\qquad
 H_\Omega=\left(|\Omega|^{-1}\int_\Omega H_j\right)_j.
\end{equation}
If $H\in W^{1,s}_0(\Omega;\lN)$, then
\begin{equation}\label{eq:Poincare-zero}
 \|H\|_{W^{1,s}(\Omega;\lN)}
 \lesssim \|\nabla H\|_{L^s(\Omega;\lN)}.
\end{equation}
If $1<t\leq s<\infty$ and
\begin{equation}\label{eq:strict-gap}
 \frac1t-\frac1s\leq\frac1d,
\end{equation}
then, for $H\in W^{1,t}(\Omega;\lN)$ and
$z\in L^t(\Omega;\lN)$,
\begin{align}
 \|H\|_{L^s(\Omega;\lN)}
 &\lesssim \|H\|_{W^{1,t}(\Omega;\lN)}, \label{eq:Sobolev}\\
 \|z\|_{W^{-1,s}_0(\Omega;\lN)}
 &\lesssim \|z\|_{L^t(\Omega;\lN)} \label{eq:negative-Sobolev}.
\end{align}
All constants are independent of $N$.
\end{lem}

\begin{proof}
We may assume that $N<\infty$, since the case $N=\infty$
follows by truncation.

We begin with the Poincar\'e estimate.
For the Rademacher variables in Lemma \ref{lem:Khintchine},
the function
$\sum_j\varepsilon_j(H_j-(H_j)_\Omega)$ has mean zero.
The scalar Poincar\'e inequality therefore gives
\[
 \left\|\sum_j\varepsilon_j(H_j-(H_j)_\Omega)
                                      \right\|_{L^s(\Omega)}
 \lesssim
 \left\|\sum_j\varepsilon_j\nabla H_j
                                \right\|_{L^s(\Omega)}.
\]
Take the $s$th power and then the expectation.
By Lemma \ref{lem:Khintchine} and Fubini's theorem,
\[
 \begin{aligned}
 \|H-H_\Omega\|_{L^s(\Omega;\lN)}^s
 &\simeq
 \E\left\|\sum_j\varepsilon_j(H_j-(H_j)_\Omega)
 \right\|_{L^s(\Omega)}^s\\
 &\lesssim
 \E\left\|\sum_j\varepsilon_j\nabla H_j
 \right\|_{L^s(\Omega)}^s
 \simeq\|\nabla H\|_{L^s(\Omega;\lN)}^s.
 \end{aligned}
\]
Since $\nabla(H-H_\Omega)=\nabla H$, this proves
\eqref{eq:Poincare}. The zero-trace estimate
\eqref{eq:Poincare-zero} follows similarly from
the scalar zero-trace Poincar\'e inequality.

Next we prove the Sobolev estimate.
Set $U=(\sum_j|H_j|^2)^{1/2}$.
The Sobolev chain rule and Cauchy--Schwarz inequality give
\[
 U\in W^{1,t}(\Omega),\qquad
 |\nabla U|\leq\left(\sum_j|\nabla H_j|^2\right)^{1/2}
 \quad\text{a.e. in }\Omega.
\]
Indeed, applying the chain rule to
$U_\eta=(\eta^2+\sum_j|H_j|^2)^{1/2}$, $\eta>0$,
we obtain, for $1\leq k\leq d$,
\[
 \partial_kU_\eta
 =\frac{\operatorname{Re}\sum_j\overline{H_j}\partial_kH_j}
        {(\eta^2+\sum_j|H_j|^2)^{1/2}},
 \qquad
 |\nabla U_\eta|
 \leq\left(\sum_j|\nabla H_j|^2\right)^{1/2}.
\]
Under \eqref{eq:strict-gap}, the scalar Sobolev embedding
$W^{1,t}(\Omega)\hookrightarrow L^s(\Omega)$ is valid.
Consequently,
\[
 \begin{aligned}
 \|H\|_{L^s(\Omega;\lN)}
 &=\|U\|_{L^s(\Omega)}\\
 &\lesssim\|U\|_{L^t(\Omega)}
                   +\|\nabla U\|_{L^t(\Omega)}\\
 &\leq\|H\|_{L^t(\Omega;\lN)}
                   +\|\nabla H\|_{L^t(\Omega;\lN)}
 =\|H\|_{W^{1,t}(\Omega;\lN)}.
 \end{aligned}
\]
This proves the positive-order estimate.
The negative-order estimate follows by duality.
\end{proof}

The following is a vector-valued version of
\cite[Lemma 5.4]{Geng-Shen-2025}.

\begin{lem}[Bogovskii and pressure]\label{lem:Bogovskii}
Let $1<s<\infty$.
On a bounded connected Lipschitz domain $\Omega$, a fixed
linear Bogovskii operator $\cB$ satisfies
\begin{equation}\label{eq:Bogovskii}
 \begin{gathered}
 \div\cB h_j=h_j,\qquad
 \cB h_j\in W^{1,s}_0(\Omega;\C^d),\qquad
 \int_\Omega h_j=0,\\
 \|(\cB h_j)_j\|_{W^{1,s}(\Omega;\lN)}
       \lesssim \|h\|_{L^s(\Omega;\lN)} .
 \end{gathered}
\end{equation}
If $\int_\Omega p_j=0$, then
\begin{equation}\label{eq:pressure-negative}
 \|p\|_{L^s(\Omega;\lN)}
 \lesssim \|(\nabla p_j)_j\|_{W^{-1,s}(\Omega;\lN(\C^d))}.
\end{equation}
\end{lem}
\begin{proof}
Choose a single Bogovskii operator on $\Omega$; see
\cite[Theorem III.3.1]{Galdi-2011}.
Apply this operator and its first derivatives to
$\sum_j\varepsilon_jh_j$.
Their scalar $L^s$ bounds, Lemma \ref{lem:Khintchine},
and Fubini's theorem give \eqref{eq:Bogovskii}.
For $h\in L^{s'}(\Omega;\lN)$, put
$v_j=\cB(h_j-(h_j)_\Omega)$. Then
\[
 \sum_j\int_\Omega p_jh_j
 =\sum_j\int_\Omega p_j\div v_j
 =-\sum_j\pair{\nabla p_j}{v_j}.
\]
Since $\|h-h_\Omega\|_{L^{s'}(\Omega;\lN)}
\leq2\|h\|_{L^{s'}(\Omega;\lN)}$, duality proves
\eqref{eq:pressure-negative}.
\end{proof}

The following is a direct consequence of
\cite[Theorem 6.1.6]{Bergh-Lofstrom-1976},
applied to diagonal operator-valued symbols
on the corresponding $\ell_N^2$ spaces.

\begin{lem}[Diagonal Mikhlin multiplier theorem]\label{lem:Mikhlin}
Let $1<q<\infty$.
Suppose that
$m_j\in C^K(\R^d\setminus\{0\};\C^{k\times m})$,
where $K=\lfloor d/2\rfloor+1$, satisfy
\begin{equation}\label{eq:Mikhlin}
 \sup_{j}\sup_{|\alpha|\le K}\sup_{\xi\neq0}
 |\xi|^{|\alpha|}\|\partial^\alpha m_j(\xi)\|
 \leq M.
\end{equation}
Then
\[
 \left\|(\cF^{-1}(m_j\widehat h_j))_j
                  \right\|_{L^q(\R^d;\ell_N^2(\C^k))}
 \lesssim_{d,q} M
 \|h\|_{L^q(\R^d;\ell_N^2(\C^m))}.
\]
\end{lem}

\begin{lem}\label{lem:negative-norm}
Let $1<q<\infty$.
For $g=(g_j)_{j=1}^N\in\dotc W^{-1,q}(\R^d;\lN)$,
\begin{equation*}
 \left\|(\nabla(-\Delta)^{-1}g_j)_{j=1}^N
            \right\|_{L^q(\R^d;\lN(\C^d))}
 \simeq_{d,q}
 \|g\|_{\dotc W^{-1,q}(\R^d;\lN)}.
\end{equation*}
\end{lem}
\begin{proof}
By Hahn--Banach there is
$F=(F_j)\in L^q(\R^d;\lN(\C^d))$ such that
\[
 g_j=-\div F_j,\qquad
 \|F\|_{L^q(\R^d;\lN(\C^d))}
     =\|g\|_{\dotc W^{-1,q}(\R^d;\lN)}.
\]
Let $R_k$ be the Riesz transform with symbol $i\xi_k/|\xi|$.
The field $G_j=\nabla(-\Delta)^{-1}g_j$ is defined by
\[
 (G_j)_k=-\sum_{\ell=1}^d R_kR_\ell(F_j)_\ell.
\]
The diagonal Riesz transform estimate, a special case of
Lemma \ref{lem:Mikhlin}, gives
\[
 \|G\|_{L^q(\R^d;\lN(\C^d))}
 \lesssim_{d,q} \|F\|_{L^q(\R^d;\lN(\C^d))}.
\]
Conversely, $g_j=-\div G_j$, so duality gives
$\|g\|_{\dotc W^{-1,q}(\R^d;\lN)}
\leq\|G\|_{L^q(\R^d;\lN(\C^d))}$.
\end{proof}

\subsection{Spaces for the resolvent problem}

In this subsection, we introduce the function spaces used in the formulation of the resolvent problem.

For pressures, taken modulo an independent constant in each
component, define
\begin{equation}\label{eq:pressure-space}
 \mathcal P_{\Lambda,\Omega}(p)=
 \inf_{p_j=p_{j,1}+p_{j,2}}
 \left\{
  \|(r_jp_{j,1})_j\|_{L^q(\Omega;\lN)}
  +\|(\nabla p_{j,2})_j\|_{L^q(\Omega;\lN)}
 \right\},
\end{equation}
where $p_{j,1}\in L^q(\Omega)$ and
$p_{j,2}\in\dotc W^{1,q}(\Omega)$.
For $h_j\in W^{-1,q}(\Omega;\C^d)$, define
\begin{equation*}
 \mathcal H_{\Lambda,\Omega}(h)=
 \inf_{h_j=F_j+\div f_j}
 \left\{
  \|F\|_{L^q(\Omega;\lN)}
  +\|(r_jf_j)_j\|_{L^q(\Omega;\lN)}
 \right\}.
\end{equation*}
The decompositions are understood in distributions, or
equivalently on $W^{1,q'}_0$ test fields. Set
\begin{align}
 \|(u,p)\|_{\cX_{\Lambda,\Omega}}
 &=\|(r_j^2u_j)_j\|_{L^q(\Omega;\lN)}
   +\|(r_j\nabla u_j)_j\|_{L^q(\Omega;\lN)}
   +\mathcal P_{\Lambda,\Omega}(p),\label{eq:Xspace}\\
 \|(h,g)\|_{\cY_{\Lambda,\Omega}}
 &=\mathcal H_{\Lambda,\Omega}(h)
   +\|(r_jg_j)_j\|_{L^q(\Omega;\lN)}
   +\|(r_j^2g_j)_j\|_{\dotc W^{-1,q}(\Omega;\lN)}.
   \label{eq:Yspace}
\end{align}
The space $\cX_{\Lambda,\Omega}$ consists of pairs with finite
norm \eqref{eq:Xspace}, with each velocity
$u_j\in W^{1,q}_0(\Omega;\C^d)$ when $\Omega$ has a boundary.
The space $\cY_{\Lambda,\Omega}$ consists of data with finite
norm \eqref{eq:Yspace}. 

\section{The whole space and half-space}\label{sec:model}

We consider the generalized system
\begin{equation}\label{eq:S-resolvent-div}
 \left\{\begin{aligned}
 -\Delta u_j+\nabla p_j+\lambda_ju_j
       &=F_j+\div f_j&&\text{in }\Omega,\\
 \div u_j&=g_j&&\text{in }\Omega,\\
 u_j&=0&&\text{on }\partial\Omega.
 \end{aligned}\right.
\end{equation}
The boundary condition is omitted when $\Omega=\R^d$.
For the data, put
\begin{equation*}
 \mathcal D_{\Lambda,\Omega}(F,f,g)
 =\|F\|_{L^q(\Omega;\lN)}
     +\|rf\|_{L^q(\Omega;\lN)}
     +\|rg\|_{L^q(\Omega;\lN)}+\|r^2g\|_{\dotc W^{-1,q}(\Omega;\lN)} .
\end{equation*}

\subsection{The whole space}

\begin{thm}\label{thm:whole}
Let $1<q<\infty$ and $\lambda_j\in\Sigma_\theta$.
If $\mathcal D_{\Lambda,\R^d}(F,f,g)<\infty$, then
\eqref{eq:S-resolvent-div} on $\R^d$ has a unique solution
$(u,p)\in\cX_{\Lambda,\R^d}$,
with pressure taken modulo constants. Moreover,
\begin{equation}\label{eq:whole-est}
 \|(u,p)\|_{\cX_{\Lambda,\R^d}}
 \lesssim_{d,q,\theta} \mathcal D_{\Lambda,\R^d}(F,f,g).
\end{equation}
\end{thm}
\begin{proof}
For each $j$, existence and uniqueness follow from
\cite[Theorem 2.1]{Geng-Shen-2025}.
The construction in its proof gives the following Fourier
representations. Let
$P(\xi)=I-\xi\otimes\xi/|\xi|^2$. Then
\begin{align*}
 \widehat u_j
 &=\frac{P(\xi)}{\lambda_j+|\xi|^2}
       (\widehat F_j+i\xi_k\widehat f_{j,\cdot k})
       -\frac{i\xi}{|\xi|^2}\widehat g_j,\\
 \widehat p_j
 &=-\frac{i\xi_\alpha}{|\xi|^2}
       (\widehat F_{j,\alpha}
                       +i\xi_k\widehat f_{j,\alpha k})
       +\widehat g_j+\frac{\lambda_j}{|\xi|^2}\widehat g_j.
\end{align*}

The sector condition implies
\begin{equation}\label{eq:sector-comparison}
 |\lambda_j+|\xi|^2|\geq c_\theta(r_j^2+|\xi|^2).
\end{equation}
The symbols
\[
 \frac{r_j^2P(\xi)}{\lambda_j+|\xi|^2},\qquad
 \frac{r_ji\xi_kP(\xi)}{\lambda_j+|\xi|^2},\qquad
 \frac{\xi_k\xi_\ell P(\xi)}{\lambda_j+|\xi|^2}
\]
satisfy \eqref{eq:Mikhlin} uniformly in $j$.
This follows by differentiating and using
\eqref{eq:sector-comparison}. Thus
Lemma \ref{lem:Mikhlin} yields
\begin{align*}
 \|r^2u\|_{L^q(\R^d;\lN)}
 &\lesssim \big(
       \|F\|_{L^q(\R^d;\lN)}
       +\|rf\|_{L^q(\R^d;\lN)}
       +\|(\nabla(-\Delta)^{-1}(r_j^2g_j))_j
                                  \|_{L^q(\R^d;\lN)}\big),\\
 \|r\nabla u\|_{L^q(\R^d;\lN)}
 &\lesssim \big(
       \|F\|_{L^q(\R^d;\lN)}
       +\|rf\|_{L^q(\R^d;\lN)}
       +\|rg\|_{L^q(\R^d;\lN)}\big).
\end{align*}
For the pressure, use the decomposition $p_j=p_{j,1}+p_{j,2}$
specified by
\begin{equation*}
 \widehat p_{j,1}
    =-\frac{i\xi_\alpha}{|\xi|^2}
                         i\xi_k\widehat f_{j,\alpha k}
                         +\widehat g_j,\quad
 \widehat p_{j,2}
    =-\frac{i\xi_\alpha}{|\xi|^2}\widehat F_{j,\alpha}
                         +\frac{\lambda_j}{|\xi|^2}\widehat g_j .
\end{equation*}
The same multiplier estimates give
\begin{align*}
 &\,\quad \|(r_jp_{j,1})_j\|_{L^q(\R^d;\lN)}
       +\|(\nabla p_{j,2})_j\|_{L^q(\R^d;\lN)}\\
 &\quad\lesssim \big(
       \|F\|_{L^q(\R^d;\lN)}
       +\|rf\|_{L^q(\R^d;\lN)}
       +\|rg\|_{L^q(\R^d;\lN)}
       +\|(\nabla(-\Delta)^{-1}(r_j^2g_j))_j\|_{L^q(\R^d;\lN)}\big).
\end{align*}
Apply Lemma \ref{lem:negative-norm} to $r^2g$ to obtain
\eqref{eq:whole-est}.
\end{proof}

\subsection{The half-space}\label{subsec:half-space}

\begin{thm}\label{thm:half}
Let $1<q<\infty$ and $\lambda_j\in\Sigma_\theta$.
If $\mathcal D_{\Lambda,\R^d_+}(F,f,g)<\infty$, then
\eqref{eq:S-resolvent-div} on $\R^d_+$ has a unique
solution $(u,p)\in\cX_{\Lambda,\R^d_+}$,
with pressure taken modulo constants. It satisfies
\begin{equation}\label{eq:half-est}
 \|(u,p)\|_{\cX_{\Lambda,\R^d_+}}
 \lesssim_{d,q,\theta}\mathcal D_{\Lambda,\R^d_+}(F,f,g).
\end{equation}
Moreover, the operator
\begin{equation*}
 S_{\Lambda,\R^d_+}(u,p)
 =\big((-\Delta u_j+\nabla p_j+\lambda_ju_j)_j,
                                      (\div u_j)_j\big)
\end{equation*}
is a bounded isomorphism from $\cX_{\Lambda,\R^d_+}$
onto $\cY_{\Lambda,\R^d_+}$, with bounds depending only on $d,q,\theta$.
\end{thm}

Following the strategy of
Geng and Shen \cite[Theorem 3.1]{Geng-Shen-2025},
we first extend the data to $\R^d$ by reflection
across $x_d=0$.
For $F_j$, we use even reflection for the tangential
components $F_{j,\alpha}$, $\alpha<d$, and odd reflection
for the normal component $F_{j,d}$.
We extend $g_j$ evenly.
For the tensor $f_j$, we use even reflection when both
indices are tangential or both equal $d$, and odd
reflection when exactly one index equals $d$.
Denote the resulting extensions by
$\widetilde F,\widetilde f,\widetilde g$.
Each $L^q(\lN)$ norm is multiplied by $2^{1/q}$.
For the homogeneous negative Sobolev norm,
even reflection also yields
\begin{equation*}
 \|r^2\widetilde g\|_{\dotc W^{-1,q}(\R^d;\lN)}
 =2^{1/q}\|r^2g\|_{\dotc W^{-1,q}(\R^d_+;\lN)},
\end{equation*}
whose proof is left to the reader.
Consequently,
\[
 \mathcal D_{\Lambda,\R^d}
       (\widetilde F,\widetilde f,\widetilde g)
 \simeq \mathcal D_{\Lambda,\R^d_+}(F,f,g).
\]
Let $(U_j,P_j)$ be the whole-space solutions given by
Theorem \ref{thm:whole}. Then
\begin{equation}\label{eq:particular-solution}
\begin{aligned}
\|(r_j^2U_j)_j\|_{L^q(\R^d;\lN)}
 +\|(r_j\nabla U_j)_j\|_{L^q(\R^d;\lN)}
 +\mathcal P_{\Lambda,\R^d}(P)\lesssim
\mathcal D_{\Lambda,\R^d_+}(F,f,g).
\end{aligned}
\end{equation}
By the reflection symmetries of the extended data
and the uniqueness of the whole-space solutions,
$U_{j,\alpha}$ is even in $x_d$ for $\alpha<d$,
whereas $U_{j,d}$ is odd.
In particular, $\Tr U_{j,d}=0$.
It remains to correct the tangential boundary trace
of $U_j$.
Set $h_j=-\Tr U'_j$ and consider the homogeneous
Stokes problem
\begin{equation}\label{eq:S-boundary-correction}
 \left\{\begin{aligned}
 -\Delta v_j+\nabla \pi_j+\lambda_jv_j
       &=0&&\text{in }\R^d_+,\\
 \div v_j&=0&&\text{in }\R^d_+,\\
 v'_j=h_j,\quad v_{j,d}&=0&&\text{on }\partial\R^d_+.
 \end{aligned}\right.
\end{equation}
Once $(v_j,\pi_j)$ is constructed, the sum
\begin{align}\label{eq:up-and-UP}
     u_j=U_j|_{\R^d_+}+v_j,\qquad
 p_j=P_j|_{\R^d_+}+\pi_j 
\end{align}
solves \eqref{eq:S-resolvent-div} with zero boundary
trace. Uniqueness follows componentwise from
\cite[Theorem 3.1]{Geng-Shen-2025}.

We now analyze the boundary correction problem
\eqref{eq:S-boundary-correction}.
Write $x=(x',y)\in\R^d_+$ and $s=|\xi'|$,
where $\xi'$ denotes the tangential Fourier variable.
For $\lambda\in\Sigma_\theta$, set
$r=|\lambda|^{1/2}$ and
$B=B_\lambda(s)=(\lambda+s^2)^{1/2}$,
with $\Rea B>0$.
Introduce
\begin{equation*}
 E_\lambda(s,y)=e^{-By},\qquad
 L_\lambda(s,y)
 =s\,\frac{e^{-By}-e^{-sy}}{B-s}.
\end{equation*}
The sector condition ensures that
\begin{equation}\label{eq:B-properties}
 |B|\simeq_\theta r+s,\qquad
 \Rea B\geq c_\theta(r+s),\qquad
 B-s=\frac{\lambda}{B+s}.
\end{equation}
For tangential boundary data
$h:\R^{d-1}\to\C^{d-1}$,
rewriting the representation formulas
in \cite[(3.10)--(3.11)]{Geng-Shen-2025} gives
\begin{align}
 \widehat v(\xi',y)
 &=\Phi_\lambda(\xi',y)\widehat h(\xi'),
       \label{eq:boundary-u}\\
 \widehat\pi(\xi',y)
 &=-i(B+s)e^{-sy}
       \left(\frac{\xi'}s\right)^{\mathsf T}\widehat h(\xi'),
       \label{eq:boundary-p}
\end{align}
where
\begin{equation}\label{eq:Phi-definition}
 \Phi_\lambda(\xi',y)=
 \begin{pmatrix}
  E_\lambda I_{d-1}
       +L_\lambda\,\dfrac{\xi'}s\otimes\dfrac{\xi'}s\\[3pt]
  iL_\lambda\left(\dfrac{\xi'}s\right)^{\mathsf T}
 \end{pmatrix}.
\end{equation}
We first establish the multiplier estimates needed
to control the boundary correction.

\begin{lem}[Derivatives of the boundary matrix]
\label{lem:boundary-symbol}
There exists $\delta>0$, depending only on $\theta$,
such that for every multi-index $\alpha$ and
$\ell=0,1,2$,
\begin{equation}\label{eq:Phi-derivatives}
 s^{|\alpha|}
 \big|\partial_{\xi'}^\alpha\partial_y^\ell
                              \Phi_\lambda(\xi',y)\big|
 \lesssim_{\theta,\alpha}
 \frac{(r+s)^\ell}{1+ry}\,e^{-\delta sy},
 \qquad s>0,\quad y>0.
\end{equation}
\end{lem}

\begin{proof}
We begin with radial derivatives.
The identities
\[
 \partial_sB=\frac sB,\qquad
 \partial_s^2B=\frac{\lambda}{B^3}
\]
and repeated differentiation yield
\begin{equation}\label{eq:B-radial-derivatives}
 |\partial_s^kB|
 \lesssim_{\theta,k}(r+s)^{1-k},\qquad
 s^k|\partial_s^kB|
 \lesssim_{\theta,k}s,\qquad k\geq1.
\end{equation}
By the chain rule,
\[
 s^k|\partial_s^kE_\lambda(s,y)|
 \lesssim_{\theta,k}
 (1+sy)^k e^{-c_\theta(r+s)y},
 \qquad k\geq0.
\]
Fix $0<\delta<c_\theta$.
Absorbing the polynomial factor into the
exponential and using
$e^{-c_\theta ry}\lesssim_\theta(1+ry)^{-1}$,
we obtain
\begin{equation}\label{eq:E-radial}
 s^k|\partial_s^kE_\lambda(s,y)|
 \lesssim_{\theta,k}
 \frac{e^{-\delta sy}}{1+ry},
 \qquad k\geq0.
\end{equation}

For $L_\lambda$, we use the integral representation
\begin{equation}\label{eq:L-integral}
 L_\lambda(s,y)
 =-sy\int_0^1
 e^{-y[(1-\tau)s+\tau B_\lambda(s)]}\,d\tau.
\end{equation}
With $c_\theta\leq1$, the phase satisfies
\[
 \Rea\big((1-\tau)s+\tau B_\lambda(s)\big)
       \geq c_\theta(s+\tau r),
 \qquad
 s^k\big|\partial_s^k[(1-\tau)s+\tau B_\lambda(s)]\big|
       \lesssim_{\theta,k} s,\quad k\geq1,
\]
uniformly for $0\leq\tau\leq1$.
Hence
\[
 s^k\left|
 \partial_s^k e^{-y[(1-\tau)s+\tau B_\lambda(s)]}
 \right|
 \lesssim_{\theta,k} (1+sy)^k e^{-c_\theta sy}e^{-c_\theta\tau ry},\quad k\geq 0.
\]
Differentiating \eqref{eq:L-integral}, we obtain
\begin{align*}
 s^k|\partial_s^kL_\lambda(s,y)|
 &\lesssim_{\theta,k}
 sy(1+sy)^k e^{-c_\theta sy}
 \int_0^1e^{-c_\theta\tau ry}\,d\tau\\
 &\lesssim_{\theta,k}
 \frac{sy(1+sy)^k}{1+ry}e^{-c_\theta sy}\lesssim_{\theta,k}
 \frac{e^{-\delta sy}}{1+ry},
\end{align*}
where we used
\[
 \int_0^1e^{-c_\theta\tau ry}\,d\tau
 \lesssim_\theta\frac1{1+ry}.
\]
For any radial function $f(s)$,
\[
 s^{|\alpha|}
 \big|\partial_{\xi'}^\alpha f(|\xi'|)\big|
 \lesssim_\alpha
 \sum_{k=0}^{|\alpha|}s^k|\partial_s^kf(s)|.
\]
Combining this with \eqref{eq:E-radial}
and the preceding estimate for $L_\lambda$ yields
\begin{equation}\label{eq:EL-tangential}
 s^{|\alpha|}
 \big(
 |\partial_{\xi'}^\alpha E_\lambda|
 +|\partial_{\xi'}^\alpha L_\lambda|
 \big)
 \lesssim_{\theta,\alpha}
 \frac{e^{-\delta sy}}{1+ry}.
\end{equation}

The normal derivatives satisfy
\begin{equation}\label{eq:EL-normal}
 \begin{aligned}
 \partial_yE_\lambda&=-BE_\lambda,&
 \partial_y^2E_\lambda&=B^2E_\lambda,\\
 \partial_yL_\lambda&=-s(L_\lambda+E_\lambda),&
 \partial_y^2L_\lambda&=s^2L_\lambda+s(s+B)E_\lambda.
 \end{aligned}
\end{equation}
By \eqref{eq:B-radial-derivatives}, the weighted tangential derivatives
of $s$ and $B$ are bounded by $C_\alpha(r+s)$,
while those of their quadratic products are
bounded by $C_\alpha(r+s)^2$.
Thus Leibniz's rule and \eqref{eq:EL-tangential}
give
\[
 s^{|\alpha|}
 \big(
 |\partial_{\xi'}^\alpha\partial_y^\ell E_\lambda|
 +|\partial_{\xi'}^\alpha\partial_y^\ell L_\lambda|
 \big)
 \lesssim_{\theta,\alpha}
 \frac{(r+s)^\ell}{1+ry}e^{-\delta sy},
 \qquad \ell=0,1,2.
\]
Finally, $\xi'/s$ and $\xi'/s\otimes\xi'/s$ are smooth
and homogeneous of degree zero away from the origin.
Their derivatives of order $|\alpha|$ are therefore
bounded by $C_\alpha s^{-|\alpha|}$.
Substituting these bounds into \eqref{eq:Phi-definition}
proves \eqref{eq:Phi-derivatives}.
\end{proof}

\begin{lem}[Four multiplier estimates]
\label{lem:four-multipliers}
Define
\begin{equation*}
 \begin{aligned}
 m_\lambda^{(1)}&=r\Phi_\lambda,&
 m_\lambda^{(2)}&=s\Phi_\lambda,\\
 m_\lambda^{(3)}&=\partial_y\Phi_\lambda,&
 m_\lambda^{(4)}&=\frac{\partial_y^2\Phi_\lambda}{r+s}.
 \end{aligned}
\end{equation*}
For $\nu=1,2,3,4$ and every multi-index $\alpha$,
\begin{equation}\label{eq:four-symbol-est}
 s^{|\alpha|}
 \big|\partial_{\xi'}^\alpha
 m_\lambda^{(\nu)}(\xi',y)\big|
 \lesssim_{\theta,\alpha}\frac1y,
 \qquad s>0,\quad y>0.
\end{equation}
The same estimate holds for
$i\xi_k\Phi_\lambda$, $1\leq k\leq d-1$,
and for the scalar symbol $se^{-sy}$.
\end{lem}

\begin{proof}
For the first three symbols,
Lemma \ref{lem:boundary-symbol} and Leibniz's rule
give the respective bounds
\[
 \frac r{1+ry}e^{-\delta sy},\qquad
 \frac s{1+ry}e^{-\delta sy},\qquad
 \frac{r+s}{1+ry}e^{-\delta sy}.
\]
For the fourth symbol, observe that
\begin{align}\label{eq:rs-1}
 \partial_s^k(r+s)^{-1}
 =\frac{(-1)^kk!}{(r+s)^{k+1}},\qquad 
 s^{|\alpha|}
 \big|\partial_{\xi'}^\alpha(r+s)^{-1}\big|
 \lesssim_\alpha\frac1{r+s}.
\end{align}
Applying Leibniz's rule and
\eqref{eq:Phi-derivatives} with $\ell=2$,
we obtain
\[
 s^{|\alpha|}
 \big|\partial_{\xi'}^\alpha
 m_\lambda^{(4)}(\xi',y)\big|
 \lesssim_{\theta,\alpha}
 \frac{r+s}{1+ry}e^{-\delta sy}.
\]
Since
\[
 \frac{r+s}{1+ry}e^{-\delta sy}
 \leq \frac1y+\frac{s e^{-\delta sy}}{1+ry}
 \lesssim_\delta\frac1y,
\]
the estimates in \eqref{eq:four-symbol-est} follow. The estimates for $i\xi_k\Phi_\lambda$ and $se^{-sy}$ follow similarly.
\end{proof}

Consequently, for each fixed $y>0$, Lemma \ref{lem:Mikhlin}
in dimension $d-1$ bounds the diagonal tangential Fourier
multipliers associated with the symbols above on
$L^q(\R^{d-1};\lN)$ by $C/y$.
These bounds are uniform in the parameter family
$(\lambda_j)_j$ and in $N$.

\begin{lem}[An integral on the half-line]\label{lem:Hardy}
Let $1<q<\infty$ and $b\geq0$. Then
\begin{equation*}
 \left\|\int_0^\infty\frac{b(t)}{y+t}\,dt
                         \right\|_{L^q_y(0,\infty)}
 \lesssim_q \|b\|_{L^q(0,\infty)}.
\end{equation*}
\end{lem}
\begin{proof}
Set $t=yz$ and apply Minkowski's inequality:
\[
 \left\|\int_0^\infty\frac{b(yz)}{1+z}\,dz
                                      \right\|_{L^q_y}
 \leq
 \left(\int_0^\infty\frac{z^{-1/q}}{1+z}\,dz\right)
                                      \|b\|_{L^q}.
\]
The last integral is finite because $1<q<\infty$.
\end{proof}

We will use these results as follows. If $T_j(\tau)$ is
a family of tangential Fourier multipliers whose symbols
satisfy the bounds $C_\alpha/\tau$ in
\eqref{eq:four-symbol-est}, then Lemma \ref{lem:Mikhlin}
and Minkowski's inequality give
\[
 \left\|\left(\int_0^\infty
             T_j(y+t)G_j(\cdot,t)\,dt\right)_j
                     \right\|_{L^q(\R^{d-1};\lN)}
 \lesssim \int_0^\infty
       \frac{\|(G_j(\cdot,t))_j\|_{L^q(\R^{d-1};\lN)}}{y+t}
       \,dt.
\]
Applying Lemma \ref{lem:Hardy} in $y$ yields
\begin{equation}\label{eq:normal-integration}
 \left\|\left(\int_0^\infty
              T_j(y+t)G_j(\cdot,t)\,dt\right)_j
                         \right\|_{L^q(\R^d_+;\lN)}
 \lesssim \|(G_j)_j\|_{L^q(\R^d_+;\lN)}.
\end{equation}
In particular, this applies to each of the symbols in
Lemma \ref{lem:four-multipliers}.

\begin{proof}[Proof of Theorem \ref{thm:half}]

The lifting identities below may first be proved for smooth
$U'_j$ with compact support in $\overline{\R^d_+}$.
The resulting estimates extend them to
$W^{1,q}(\R^d_+)$ by density.

\emph{The velocity estimates.}
Since $h_j=-\Tr U'_j$, \eqref{eq:boundary-u} and integration
of the derivative of
$\Phi_{\lambda_j}(\xi',y+t)\widehat U'_j(\xi',t)$
with respect to $t$ give
\begin{equation}\label{eq:trace-extension}
 \widehat v_j(\xi',y)
 =\int_0^\infty
 \left\{
  (\partial_y\Phi_{\lambda_j})(\xi',y+t)\widehat U'_j(\xi',t)
  +\Phi_{\lambda_j}(\xi',y+t)\partial_t\widehat U'_j(\xi',t)
 \right\}\,dt.
\end{equation}
Multiplying \eqref{eq:trace-extension} by $r_j^2$ yields
\begin{equation}\label{eq:velocity-interior}
 r_j^2\widehat v_j
 =\int_0^\infty
 \left\{
  (\partial_y\Phi_{\lambda_j})\,r_j^2\widehat U'_j
  +(r_j\Phi_{\lambda_j})\,r_j\partial_t\widehat U'_j
 \right\}\,dt.
\end{equation}
Here and in the two derivative formulas below, the matrices are
evaluated at $(\xi',y+t)$ and the functions involving
$U'_j$ at $(\xi',t)$.
The two matrices in \eqref{eq:velocity-interior} are
$m_{\lambda_j}^{(3)}$ and $m_{\lambda_j}^{(1)}$.
Thus \eqref{eq:normal-integration} bounds
$\|(r_j^2v_j)_j\|_{L^q(\R^d_+;\lN)}$ by
\[
 C\big(
 \|(r_j^2U'_j)_j\|_{L^q(\R^d_+;\lN)}
 +\|(r_j\partial_yU'_j)_j\|_{L^q(\R^d_+;\lN)}\big).
\]

For $1\leq k\leq d-1$, multiplication of
\eqref{eq:trace-extension} by $r_ji\xi_k$ gives
\begin{equation*}
 r_j\widehat{\partial_kv_j}
 =\int_0^\infty
 \left\{
  (\partial_y\Phi_{\lambda_j})\,r_j\widehat{\partial_kU'_j}
  +(i\xi_k\Phi_{\lambda_j})\,r_j\partial_t\widehat U'_j
 \right\}\,dt.
\end{equation*}
Both symbols satisfy the bounds in
Lemma \ref{lem:four-multipliers}, so
\eqref{eq:normal-integration} bounds this derivative by
the $L^q(\R^d_+;\lN)$ norms of
$(r_j\partial_kU'_j)_j$ and $(r_j\partial_yU'_j)_j$.

For the normal derivative, differentiate
\eqref{eq:trace-extension} in $y$ to write
\begin{equation}\label{eq:normal-interior}
 r_j\widehat{\partial_yv_j}
 =\int_0^\infty
 \left\{
  \frac{\partial_y^2\Phi_{\lambda_j}}{r_j+s}
           (r_j^2\widehat U'_j+r_js\widehat U'_j)
  +(\partial_y\Phi_{\lambda_j})\,r_j\partial_t\widehat U'_j
 \right\}\,dt.
\end{equation}
The matrices are now $m_{\lambda_j}^{(4)}$ and
$m_{\lambda_j}^{(3)}$.
The input $r_js\widehat U'_j$ is the tangential Fourier
transform of $r_j|D'|U'_j$, where $|D'|$ has symbol
$|\xi'|=s$.
If $R_k$ denotes the tangential Riesz transform with
symbol $i\xi_k/s$, then
$|D'|=-\sum_{k=1}^{d-1}R_k\partial_k$.
Lemma \ref{lem:Mikhlin} therefore gives
\begin{equation}\label{eq:Ga-est}
 \|(r_j|D'|U'_j)_j\|_{L^q(\R^d_+;\lN)}
 \lesssim \|(r_j\nabla'U'_j)_j\|_{L^q(\R^d_+;\lN)}.
\end{equation}
Apply \eqref{eq:normal-integration} to
\eqref{eq:normal-interior} and then use
\eqref{eq:Ga-est}.
Combining the velocity and derivative estimates, we obtain
\begin{equation}\label{eq:boundary-velocity-est}
 \begin{split}
 &\|(r_j^2v_j)_j\|_{L^q(\R^d_+;\lN)}
       +\|(r_j\nabla v_j)_j\|_{L^q(\R^d_+;\lN)}\\
 &\qquad\lesssim 
       \|(r_j^2U'_j)_j\|_{L^q(\R^d_+;\lN)}
       +\|(r_j\nabla U'_j)_j\|_{L^q(\R^d_+;\lN)}.
 \end{split}
\end{equation}

\emph{The pressure estimate.}
We use the decomposition
\begin{equation*}
 \widehat\pi_{j,1}
       =\frac{s}{r_j+s}\widehat\pi_j,\qquad
 \widehat\pi_{j,2}
       =\frac{r_j}{r_j+s}\widehat\pi_j.
\end{equation*}
Let $z_j$ be the Poisson extension of $h_j$:
$\widehat z_j(\xi',y)=e^{-sy}\widehat h_j(\xi')$.
The pressure formula \eqref{eq:boundary-p} gives
\begin{equation*}
 \widehat{r_j\pi_{j,1}}
       =c_j(\xi')\,r_js\widehat z_j,\qquad
 \widehat{\nabla\pi_{j,2}}
       =\begin{pmatrix}i\xi'/s\\-1\end{pmatrix}
                      c_j(\xi')\,r_js\widehat z_j,
\end{equation*}
where the row symbol is
\[
 c_j(\xi')=-i\,\frac{B_{\lambda_j}(s)+s}{r_j+s}
                           \left(\frac{\xi'}s\right)^{\mathsf T}.
\]
The bounds in \eqref{eq:B-radial-derivatives} and
\eqref{eq:rs-1} show that
\begin{equation*}
 s^{|\alpha|}|\partial_{\xi'}^\alpha c_j(\xi')|
 \lesssim_{\theta,\alpha} 1.
\end{equation*}
The column $(i\xi'/s,-1)^{\mathsf T}$ satisfies the same
zero-order bounds.
Hence Lemma \ref{lem:Mikhlin}, applied at each height,
implies
\begin{equation}\label{eq:pressure-to-Poisson}
 \|(r_j\pi_{j,1})_j\|_{L^q(\R^d_+;\lN)}
 +\|(\nabla\pi_{j,2})_j\|_{L^q(\R^d_+;\lN)}
 \lesssim \|(r_j|D'|z_j)_j\|_{L^q(\R^d_+;\lN)}.
\end{equation}
Integrating the derivative of
$e^{-s(y+t)}\widehat U'_j(\xi',t)$ gives
\begin{equation}\label{eq:Poisson-interior}
 r_js\widehat z_j(\xi',y)
 =\int_0^\infty se^{-s(y+t)}
       \big(-r_js\widehat U'_j(\xi',t)
                    +r_j\partial_t\widehat U'_j(\xi',t)\big)\,dt.
\end{equation}
The symbol $se^{-s\tau}$ satisfies the bounds
$C_\alpha/\tau$ in Lemma \ref{lem:four-multipliers}.
Thus \eqref{eq:normal-integration} and \eqref{eq:Ga-est}
yield
\begin{equation}\label{eq:Poisson-extension-est}
 \begin{split}
 \|(r_j|D'|z_j)_j\|_{L^q(\R^d_+;\lN)}
 &\lesssim \big(
       \|(r_j|D'|U'_j)_j\|_{L^q(\R^d_+;\lN)}
       +\|(r_j\partial_yU'_j)_j\|_{L^q(\R^d_+;\lN)}\big)\\
 &\lesssim \|(r_j\nabla U'_j)_j\|_{L^q(\R^d_+;\lN)}.
 \end{split}
\end{equation}
Combining this with \eqref{eq:pressure-to-Poisson}, we obtain
\begin{equation*}
 \mathcal P_{\Lambda,\R^d_+}(\pi)
 \lesssim \|(r_j\nabla U'_j)_j\|_{L^q(\R^d_+;\lN)}.
\end{equation*}
Finally, combining the above estimate with
\eqref{eq:boundary-velocity-est}, \eqref{eq:up-and-UP} and
\eqref{eq:particular-solution} yields \eqref{eq:half-est}.

\emph{The isomorphism.}
For $h=F+\div f$, the preceding construction gives a solution
satisfying \eqref{eq:half-est}.
Uniqueness makes the solution independent of the
decomposition.
Taking the infimum over all decompositions proves
surjectivity and the bound for $S_{\Lambda,\R^d_+}^{-1}$.

To prove boundedness of $S_{\Lambda,\R^d_+}$, write
$p_j=p_{j,1}+p_{j,2}$ and choose
\[
 F_j=\lambda_ju_j+\nabla p_{j,2},\qquad
 f_j=-\nabla u_j+p_{j,1}I.
\]
By the definition of $\mathcal H_{\Lambda,\R^d_+}$,
\[
 \begin{split}
 \mathcal H_{\Lambda,\R^d_+}
       ((-\Delta u_j+\nabla p_j+\lambda_ju_j)_j)&\lesssim \big(
       \|(r_j^2u_j)_j\|_{L^q(\R^d_+;\lN)}
       +\|(r_j\nabla u_j)_j\|_{L^q(\R^d_+;\lN)}\\
 &\,\quad 
       +\|(\nabla p_{j,2})_j\|_{L^q(\R^d_+;\lN)}
       +\|(r_jp_{j,1})_j\|_{L^q(\R^d_+;\lN)}\big).
 \end{split}
\]
For the divergence datum,
\[
 \|(r_j\div u_j)_j\|_{L^q(\R^d_+;\lN)}
 \lesssim \|(r_j\nabla u_j)_j\|_{L^q(\R^d_+;\lN)}.
\]
Since each $u_j$ has zero trace, integration by parts
against tests of arbitrary trace gives
\[
 \sum_j\pair{r_j^2\div u_j}{\varphi_j}
       =-\sum_j\int_{\R^d_+}r_j^2u_j\cdot\nabla\varphi_j.
\]
Consequently,
\[
 \|(r_j^2\div u_j)_j\|_{\dotc W^{-1,q}(\R^d_+;\lN)}
 \leq\|(r_j^2u_j)_j\|_{L^q(\R^d_+;\lN)}.
\]
Taking the infimum over the pressure decompositions
proves the forward operator bound.
\end{proof}

\section{The region above a Lipschitz graph}\label{sec:graph}

Let
\[
 \mathbb H_\psi=\{x=(x',x_d):x_d>\psi(x')\},
 \qquad \psi:\R^{d-1}\to\R,\qquad a=\nabla'\psi .
\]
Only $a\in L^\infty$ will be used. Set
\[
 T(y',y_d)=(y',y_d+\psi(y')),\qquad \Psi=T^{-1}.
\]
Both maps are Lipschitz, and $\det DT=1$. Moreover,
\[
 DT=
 \begin{pmatrix}
 I_{d-1}&0\\
 a^{\mathsf T}&1
 \end{pmatrix},
 \qquad
 \big((DT)^{-1}\big)^{\mathsf T}=
 \begin{pmatrix}
 I_{d-1}&-a\\
 0&1
 \end{pmatrix}.
\]
For fields on $\mathbb H_\psi$ write
$v_j=u_j\circ T$, $\pi_j=p_j\circ T$,
$\widehat F_j=F_j\circ T$, $\widehat g_j=g_j\circ T$,
and $\widetilde f_j=f_j\circ T$. Define the transformed flux by
\begin{equation}\label{eq:flux-transform}
 \widehat f_{j,\alpha k}=\widetilde f_{j,\alpha k}\quad(k<d),
 \qquad
 \widehat f_{j,\alpha d}=
   \widetilde f_{j,\alpha d}
       -\sum_{k=1}^{d-1}a_k\widetilde f_{j,\alpha k}.
\end{equation}
In matrix notation,
$\widehat f=\widetilde f\big((DT)^{-1}\big)^{\mathsf T}$.
This is the flux transformation associated with the change of
coordinates: for every zero-trace test $\varphi$,
\[
 -\int_{\mathbb H_\psi}f_{\alpha k}\partial_{x_k}\varphi_\alpha
 =-\int_{\half}\widehat f_{\alpha k}
                     \partial_{y_k}(\varphi_\alpha\circ T).
\]
In particular, $(\div_xf)\circ T=\div_y\widehat f$ in
distributions.

\begin{lem}[Change of coordinates]\label{lem:coordinates}
Suppose $\|a\|_\infty\leq1$. Then
\begin{align}
 \|(v,\pi)\|_{\cX_{\Lambda,\half}}
 &\simeq_{d,q}\|(u,p)\|_{\cX_{\Lambda,\mathbb H_\psi}},\notag
\\
 \|(\widehat h,\widehat g)\|_{\cY_{\Lambda,\half}}
 &\simeq_{d,q}\|(h,g)\|_{\cY_{\Lambda,\mathbb H_\psi}},
 \qquad \widehat h=\widehat F+\div\widehat f,\quad
 h=F+\div f .\label{eq:Y-coordinates}
\end{align}
\end{lem}
\begin{proof}
Since $\det DT=1$, all $L^q(\lN)$ norms are preserved
under composition with $T$. The chain rule gives the
equivalence of gradient norms, while duality yields
\[
 \|r^2\widehat g\|_{\dotc W^{-1,q}(\half;\lN)}
 \simeq
 \|r^2g\|_{\dotc W^{-1,q}(\mathbb H_\psi;\lN)}.
\]
For the $\mathcal H$ norm, each decomposition
$h=F+\div f$ transforms into
$\widehat h=\widehat F+\div\widehat f$ by
\eqref{eq:flux-transform}, with comparable weighted
data norms. This completes the proof.
\end{proof}

Under the change of coordinates, system
\eqref{eq:S-resolvent-div} on $\mathbb H_\psi$
becomes the following system on $\half$:
\begin{equation*}
 \left\{\begin{aligned}
 -\Delta v_j+\nabla\pi_j+\lambda_jv_j
   &=\widehat F_j+\div\widehat f_j-R_m(v_j,\pi_j),\\
 \div v_j&=\widehat g_j-R_g(v_j),\\
 v_j&=0\quad\text{on }\partial\half,
 \end{aligned}\right.
\end{equation*}
where we set $a_d=0$ and sum over $k=1,\ldots,d-1$:
\begin{align}
 [R_m(v_j,\pi_j)]_\alpha
 &=\partial_d(a_k\partial_kv_{j,\alpha})
   +\partial_k(a_k\partial_dv_{j,\alpha})
   -\partial_d(|a|^2\partial_dv_{j,\alpha})
   -\partial_d(a_\alpha\pi_j),\label{eq:Rm}\\
 R_g(v_j)&=-\partial_d(a_kv_{j,k}).\notag
\end{align}

\begin{lem}[Perturbation estimate]\label{lem:perturbation}
If $\|a\|_\infty\leq1$, then
\begin{equation*}
 \|(R_m(v,\pi),R_g(v))\|_{\cY_{\Lambda,\half}}
 \lesssim \|a\|_\infty\|(v,\pi)\|_{\cX_{\Lambda,\half}}.
\end{equation*}
\end{lem}
\begin{proof}
The first three terms in \eqref{eq:Rm} are divergences of
matrices whose $r_j$-weighted norms are bounded by
$C(\|a\|_\infty+\|a\|_\infty^2)\|r\nabla v\|_{L^q(\R^d_+;\lN)}$.
For pressure, choose $\pi_j=\pi_{j,1}+\pi_{j,2}$.
The term $-\partial_d(a_\alpha\pi_{j,1})$ is a
divergence-form force with weighted norm at most
$\|a\|_\infty\|r\pi_1\|_{L^q(\R^d_+;\lN)}$.
Since $a$ is independent of $y_d$, the remaining term is
the $L^q$ force $-a_\alpha\partial_d\pi_{j,2}$, bounded
by $\|a\|_\infty\|\nabla\pi_2\|_{L^q(\R^d_+;\lN)}$.
Taking the infimum over the pressure decomposition proves the
bound for $R_m$.

For divergence,
\[
 \begin{aligned}
 \|rR_g(v)\|_{L^q(\R^d_+;\lN)}
 &\leq\|a\|_\infty\|r\nabla v\|_{L^q(\R^d_+;\lN)},\\
 \|r^2R_g(v)\|_{\dotc W^{-1,q}(\R^d_+;\lN)}
 &\leq\|a\|_\infty\|r^2v\|_{L^q(\R^d_+;\lN)}.
 \end{aligned}
\]
The second bound follows by pairing
$-\partial_d(a_k r_j^2v_{j,k})$ with homogeneous tests.
There is no boundary contribution because $v$ has zero
trace. Combining these estimates proves the lemma.
\end{proof}

\begin{thm}[Small Lipschitz graph]\label{thm:graph}
Let $1<q<\infty$ and $\lambda_j\in\Sigma_\theta$.
There exists $\varepsilon=\varepsilon(d,q,\theta)>0$
such that, if $\|a\|_\infty\leq\varepsilon$,
then $S_{\Lambda,\mathbb H_\psi}$ is a bounded
isomorphism from $\cX_{\Lambda,\mathbb H_\psi}$
onto $\cY_{\Lambda,\mathbb H_\psi}$.
In particular, if
$\mathcal D_{\Lambda,\mathbb H_\psi}(F,f,g)<\infty$,
then \eqref{eq:S-resolvent-div} on $\mathbb H_\psi$
has a unique solution
$(u,p)\in\cX_{\Lambda,\mathbb H_\psi}$,
with pressure taken modulo constants, and
\begin{equation*}
 \|(u,p)\|_{\cX_{\Lambda,\mathbb H_\psi}}
 \lesssim_{d,q,\theta}
 \mathcal D_{\Lambda,\mathbb H_\psi}(F,f,g).
\end{equation*}
\end{thm}
\begin{proof}
Under the change of coordinates, the system on
$\mathbb H_\psi$ is equivalent to
\[
 (S_{\Lambda,\half}+R)(v,\pi)
 =(\widehat h,\widehat g),
 \qquad R(v,\pi)=(R_m(v,\pi),R_g(v)).
\]
By Theorem \ref{thm:half} and Lemma \ref{lem:perturbation},
\[
 \|S_{\Lambda,\half}^{-1}R\|_{\cL(\cX_{\Lambda,\half})}
 \lesssim_{d,q,\theta}\|a\|_\infty.
\]
Choose $\varepsilon=\varepsilon(d,q,\theta)>0$ so that
the left-hand side is at most $1/2$.
Then
$S_{\Lambda,\half}+R$ is invertible, with inverse given
by the Neumann series
\[
 (S_{\Lambda,\half}+R)^{-1}
 =\sum_{\nu=0}^\infty
 (-S_{\Lambda,\half}^{-1}R)^\nu
 S_{\Lambda,\half}^{-1}.
\]
The series converges in operator norm, with a uniform
bound. The conclusion follows from Lemma
\ref{lem:coordinates}.
\end{proof}

\section{A bounded \texorpdfstring{$C^1$}{C1} domain}\label{sec:bounded}

Throughout this section $\Omega$ is a bounded connected
$C^1$ domain in $\R^d$.
For every $\delta>0$, there exists $r_0>0$ such that,
for each $z\in\partial\Omega$, there are a graph domain
$D_z$ and a rigid motion $R_z$ taking $z$ to the origin,
with
\begin{equation}\label{eq:local-graph-agreement}
 \Omega\cap B(z,2r_0)=D_z\cap B(z,2r_0),\qquad
 \partial\Omega\cap B(z,2r_0)=\partial D_z\cap B(z,2r_0),
\end{equation}
and
\[
 R_z(D_z)=\mathbb H_{\psi_z},\qquad
 \psi_z\in C^1(\R^{d-1}),\qquad
 \psi_z(0)=0,\quad \nabla'\psi_z(0)=0,\quad
 \|\nabla'\psi_z\|_\infty\leq\delta.
\]
This is the local graph description used in
\cite[Section 5]{Geng-Shen-2025}.

Choose finitely many boundary balls $B(z_\ell,r_0)$ and
interior balls so that the smaller balls cover
$\overline\Omega$ and each doubled interior ball lies in
$\Omega$.
For a boundary ball, set $D_\ell=D_{z_\ell}$ and choose
$\phi_\ell\in C^\infty_c(B(z_\ell,2r_0))$ with
$\phi_\ell=1$ on $B(z_\ell,r_0)$.
For an interior ball, set $D_\ell=\R^d$ and choose the
cutoff in the same way.
Thus each localized field $\phi_\ell u$ is regarded
as a function on $D_\ell$, extended by zero outside
$B(z_\ell,2r_0)$.
The radii and cutoff bounds depend only on $\delta$ and
the $C^1$ character of $\Omega$.

The goal of this section is to prove the following theorem, which contains Theorem \ref{thm:main} as a special case with $f=0$ and $g=0$.

\begin{thm}[Generalized bounded-domain estimate]\label{thm:general}
Let $\Omega\subset\R^d$ be a bounded connected $C^1$ domain,
$1<q<\infty$ and $\lambda_j\in\Sigma_\theta$.
For $F_j\in L^q(\Omega;\C^d)$,
$f_j\in L^q(\Omega;\C^{d\times d})$ and
$g_j\in L^q_0(\Omega)$, problem \eqref{eq:S-resolvent-div}
has a unique solution
$u_j\in W^{1,q}_0(\Omega;\C^d)$,
$p_j\in L^q_0(\Omega)$. With
$\rho_j=(1+|\lambda_j|)^{1/2}$, it satisfies
\begin{equation}\label{eq:general-est}
 \|\rho\nabla u,\rho^2u,p\|_{L^q(\Omega;\lN)}
 \lesssim_{d,q,\theta,\Omega}
 \|F,\rho f,\rho^2g\|_{L^q(\Omega;\lN)}.
\end{equation}
The conclusion also holds for $\lambda_j=0$,
with $\rho_j=1$.
\end{thm}

\subsection{The localized divergence and pressure}

The following lemma is a vector-valued version of
\cite[Lemma 5.2 and Remark 5.3]{Geng-Shen-2025}.

\begin{lem}[Localized divergence]\label{lem:local-div}
Suppose $u_j\in W^{1,s}_0(\Omega;\C^d)$ and $\div u_j=0$.
For any cutoff $\phi$ and its associated model domain $D$,
\begin{equation}\label{eq:local-div}
 \|(\div(\phi u_j))_j\|_{\dotc W^{-1,s}(D;\lN)}
 \lesssim_\phi
 \|(u_j)_j\|_{W^{-1,s}_0(\Omega;\lN)}.
\end{equation}
\end{lem}
\begin{proof}
Let $H=(H_j)\in\dotc W^{1,s'}(D;\lN)$.
Choose a bounded connected Lipschitz patch $K$ in the
common coordinate neighborhood containing
$\supp\nabla\phi\cap\Omega$, and set
$H_K=(|K|^{-1}\int_K H_j)_j$.
Since $\div u_j=0$ and $\phi u_j$ has zero trace,
\[
 \begin{aligned}
 \left|\sum_j\pair{\div(\phi u_j)}{H_j}\right|
 &=\left|\sum_j\int_\Omega
 u_j\cdot\nabla\phi\,(H_j-(H_K)_j)\right|\\
 &\leq
 \|(u_j)_j\|_{W^{-1,s}_0(\Omega;\lN)}
 \|\nabla\phi(H-H_K)\|_{W^{1,s'}(\Omega;\lN)}.
 \end{aligned}
\]
Extending $\nabla\phi(H-H_K)$ by zero outside the patch,
the product rule and \eqref{eq:Poincare} give
\[
 \begin{aligned}
 \|\nabla\phi(H-H_K)\|_{W^{1,s'}(\Omega;\lN)}
 &\lesssim_\phi
 \|H-H_K\|_{W^{1,s'}(K;\lN)}\\
 &\lesssim_\phi
 \|\nabla H\|_{L^{s'}(K;\lN)}\leq\|\nabla H\|_{L^{s'}(D;\lN)}.
 \end{aligned}
\]
The desired estimate follows by duality.
\end{proof}

By the equation
$\nabla p=\Delta u-\lambda u+F+\div f$
and \eqref{eq:pressure-negative},
\begin{equation}\label{eq:bounded-pressure}
 \|p\|_{L^s(\Omega;\lN)}
 \lesssim
 \|F,f,\nabla u\|_{L^s(\Omega;\lN)}
 +\|r^2u\|_{W^{-1,s}_0(\Omega;\lN)}.
\end{equation}
Indeed,
\[
 \begin{aligned}
 \|\Delta u\|_{W^{-1,s}(\Omega;\lN)}
 &\leq \|\nabla u\|_{L^s(\Omega;\lN)},\\
 \|F\|_{W^{-1,s}(\Omega;\lN)}
 &\lesssim \|F\|_{L^s(\Omega;\lN)},\\
 \|\div f\|_{W^{-1,s}(\Omega;\lN)}
 &\leq \|f\|_{L^s(\Omega;\lN)},\\
 \|\lambda u\|_{W^{-1,s}(\Omega;\lN)}
 &\leq \|r^2u\|_{W^{-1,s}_0(\Omega;\lN)}.
 \end{aligned}
\]
The first and third inequalities follow by duality,
the second from the embedding $L^s\hookrightarrow W^{-1,s}$,
and the last from the embedding
$W^{-1,s}_0(\Omega;\lN)\hookrightarrow W^{-1,s}(\Omega;\lN)$.

\begin{lem}[Large parameters with a negative-norm remainder]
\label{lem:local-est}
Let $1<s<\infty$ and $\lambda_j\in\Sigma_\theta$.
There exists $M=M(d,s,\theta,\Omega)>1$ such that,
if $|\lambda_j|\geq M$ for all $j$, then every solution
of \eqref{eq:S-resolvent-div} with $g=0$,
$u_j\in W^{1,s}_0(\Omega)$ and $p_j\in L^s_0(\Omega)$
satisfies
\begin{equation}\label{eq:large-remainder}
 \|r\nabla u,r^2u\|_{L^s(\Omega;\lN)}
 \lesssim_{d,s,\theta,\Omega}
 \|F,rf\|_{L^s(\Omega;\lN)}
 +\|r^2u\|_{W^{-1,s}_0(\Omega;\lN)}.
\end{equation}
\end{lem}
\begin{proof}
Choose the boundary charts sufficiently flat to apply
Theorem \ref{thm:graph} at exponent $s$.
For each cutoff $\phi$ and its associated model domain $D$,
the localized fields satisfy
\begin{equation*}
 \left\{\begin{aligned}
 -\Delta(\phi u_j)+\nabla(\phi p_j)+\lambda_j\phi u_j
      &=\widetilde F_j+\div\widetilde f_j
        &&\text{in }D,\\
 \div(\phi u_j)&=u_j\cdot\nabla\phi
        &&\text{in }D,\\
 \phi u_j&=0&&\text{on }\partial D,
 \end{aligned}\right.
\end{equation*}
where
\[
 \begin{split}
 \widetilde F_{j,\alpha}
 &=\phi F_{j,\alpha}-f_{j,\alpha k}\partial_k\phi
   +p_j\partial_\alpha\phi+u_{j,\alpha}\Delta\phi,\\
 \widetilde f_{j,\alpha k}
 &=\phi f_{j,\alpha k}-2u_{j,\alpha}\partial_k\phi.
 \end{split}
\]
Apply Theorem \ref{thm:whole} or \ref{thm:graph} on $D$.
The $r$-weighted $L^s$ norm of the divergence datum
is bounded by $C_\phi\|ru\|_{L^s(\Omega;\lN)}$,
while its $r^2$-weighted homogeneous negative norm
is controlled by Lemma \ref{lem:local-div}.
Summing over the finite cover and recovering $\nabla u,u$
on the smaller sets where $\phi=1$, we obtain
\[
 \begin{split}
 \|r\nabla u,r^2u\|_{L^s(\Omega;\lN)}
 \leq C\big(&\|F,rf,p,ru\|_{L^s(\Omega;\lN)}+\|r^2u\|_{W^{-1,s}_0(\Omega;\lN)}\big).
 \end{split}
\]
By \eqref{eq:bounded-pressure} and $r_j\geq1$,
\[
 \begin{split}
 \|r\nabla u,r^2u\|_{L^s(\Omega;\lN)}
 \leq C\big(&\|F,rf,\nabla u,ru\|_{L^s(\Omega;\lN)}+\|r^2u\|_{W^{-1,s}_0(\Omega;\lN)}\big).
 \end{split}
\]
Since $r_j\geq\sqrt{M}$, we have
\[
 \|\nabla u,ru\|_{L^s(\Omega;\lN)}
 \leq M^{-1/2}
 \|r\nabla u,r^2u\|_{L^s(\Omega;\lN)}.
\]
Choosing $M$ sufficiently large to absorb this term
into the left-hand side gives \eqref{eq:large-remainder}.
\end{proof}

\subsection{Energy and the finite bootstrap}

The following $L^2$ estimate is standard; we include
a proof for completeness.

\begin{lem}[The $L^2$ estimate]\label{lem:energy}
Let $\lambda_j\in\Sigma_\theta\cup\{0\}$.
Every solution of \eqref{eq:S-resolvent-div} with $g=0$
satisfies
\begin{equation*}
 \|\rho\nabla u,\rho^2u\|_{L^2(\Omega;\lN)}
 \lesssim_{d,\theta,\Omega}
 \|F,\rho f\|_{L^2(\Omega;\lN)},
\end{equation*}
where $\rho_j=(1+|\lambda_j|)^{1/2}$.
\end{lem}

\begin{proof}
For each $j$, test the equation by $\overline u_j$.
Let $\alpha_j=\arg\lambda_j$ if $\lambda_j\neq0$,
and $\alpha_j=0$ otherwise.
Multiplying the resulting identity by
$e^{-i\alpha_j/2}$ and taking real parts gives
\[
 \sin(\theta/2)
 \left(\|\nabla u_j\|_{L^2(\Omega)}^2
       +|\lambda_j|\|u_j\|_{L^2(\Omega)}^2\right)
 \leq
 \left|\int_\Omega F_j\cdot\overline u_j
                  -f_j:\nabla\overline u_j\right|.
\]
Since $\rho_j^2=1+|\lambda_j|$, Poincar\'e's inequality
and H\"older's inequality yield
\[
 \begin{aligned}
 &\|\nabla u_j\|_{L^2(\Omega)}^2
   +\rho_j^2\|u_j\|_{L^2(\Omega)}^2\\
 &\quad\lesssim_{\theta,\Omega}
   \|F_j\|_{L^2(\Omega)}\|u_j\|_{L^2(\Omega)}
   +\|f_j\|_{L^2(\Omega)}\|\nabla u_j\|_{L^2(\Omega)}\\
 &\quad\leq
   \left(\rho_j^{-1}\|F_j\|_{L^2(\Omega)}
         +\|f_j\|_{L^2(\Omega)}\right)
   \left(\|\nabla u_j\|_{L^2(\Omega)}^2
         +\rho_j^2\|u_j\|_{L^2(\Omega)}^2\right)^{1/2}.
 \end{aligned}
\]
Canceling the common square-root factor, multiplying by $\rho_j$, and summing the squares over $j$
completes the proof.
\end{proof}

We next apply the finite Sobolev bootstrap of
\cite[Lemma 5.6]{Geng-Shen-2025} in the vector-valued setting.

\begin{prop}[A priori estimate for $q\geq2$]\label{prop:q-large}
Let $2\leq q<\infty$ and $\lambda_j\in\Sigma_\theta\cup\{0\}$.
Every solution of \eqref{eq:S-resolvent-div} with $g=0$
satisfies
\begin{equation}\label{eq:solenoidal-est}
 \|\rho\nabla u,\rho^2u\|_{L^q(\Omega;\lN)}
 \lesssim_{d,q,\theta,\Omega}
 \|F,\rho f\|_{L^q(\Omega;\lN)},
\end{equation}
where $\rho_j=(1+|\lambda_j|)^{1/2}$.
\end{prop}

\begin{proof}
The case $q=2$ follows from Lemma \ref{lem:energy}.
For $q>2$, choose a finite sequence
\begin{equation}\label{eq:exponent-chain}
 2=s_0<s_1<\cdots<s_K=q,\qquad
 \frac1{s_{k-1}}-\frac1{s_k}\leq \frac1d .
\end{equation}
Choose $M>\max_{1\leq k\leq K}M_{s_k}$, where
$M_{s_k}$ are given by Lemma \ref{lem:local-est}.
We consider separately the indices with
$|\lambda_j|\geq M$ and $|\lambda_j|<M$.

\emph{Large parameters.}
For $|\lambda_j|\geq M$, we have $\rho_j\simeq r_j$.
By Lemma \ref{lem:local-est}, \eqref{eq:negative-Sobolev} and
\eqref{eq:exponent-chain},
\[
 \begin{split}
 \|r\nabla u,r^2u\|_{L^{s_k}(\Omega;\lN)}
 \lesssim
 &\|F,rf\|_{L^{s_k}(\Omega;\lN)}+\|r^2u\|_{L^{s_{k-1}}(\Omega;\lN)}.
 \end{split}
\]
Starting from the $L^2$ estimate in Lemma
\ref{lem:energy} and iterating over $k=1,\ldots,K$,
we obtain \eqref{eq:solenoidal-est} for the large
parameters, since $\Omega$ is bounded.

\emph{Bounded parameters.}
For $|\lambda_j|<M$, set $\mu_j=\lambda_j+2M$.
Then $\mu_j\in\Sigma_\theta$ and
$M<|\mu_j|<3M$. Rewrite the equation as
\[
 -\Delta u_j+\nabla p_j+\mu_ju_j
 =F_j+\div f_j+2Mu_j.
\]
Applying the large-parameter estimate at exponent
$s_k$ to this equation gives
\begin{equation}\label{eq:shift-est}
 \|\nabla u,u\|_{L^{s_k}(\Omega;\lN)}
 \lesssim_M
 \|F,f,u\|_{L^{s_k}(\Omega;\lN)}.
\end{equation}
By \eqref{eq:Poincare-zero} and
\eqref{eq:Sobolev},
\[
 \|u\|_{L^{s_k}(\Omega;\lN)}
 \lesssim
 \|\nabla u\|_{L^{s_{k-1}}(\Omega;\lN)}.
\]
Combining this with \eqref{eq:shift-est} and
starting from Lemma \ref{lem:energy}, we obtain
the desired estimate by induction on $k$.
Since $1\leq\rho_j\leq(1+M)^{1/2}$, this proves
\eqref{eq:solenoidal-est} for the bounded parameters.

Combining the two cases completes the proof.
\end{proof}

\subsection{Solvability and the range \texorpdfstring{$1<q<2$}{1 < q < 2}}

\begin{proof}[Proof of Theorem \ref{thm:general}]
By \cite[Theorem 5.1 and Remark 5.7]{Geng-Shen-2025},
the existence and uniqueness of $(u_j,p_j)$ are known
for each $j$, including when $\lambda_j=0$.
It remains to prove the uniform $\ell^2$ estimates.

\emph{The case $g=0$ and $q\geq2$.}
Proposition \ref{prop:q-large} gives
\[
 \|\rho\nabla u,\rho^2u\|_{L^q(\Omega;\lN)}
 \lesssim_{d,q,\theta,\Omega}
 \|F,\rho f\|_{L^q(\Omega;\lN)}.
\]
By \eqref{eq:bounded-pressure} and the embedding
$L^q\hookrightarrow W^{-1,q}_0$, we also have
\begin{equation}\label{eq:pressure-global}
 \begin{aligned}
 \|p\|_{L^q(\Omega;\lN)}
 &\lesssim
 \|F,f,\nabla u\|_{L^q(\Omega;\lN)}
 +\|\rho^2u\|_{W^{-1,q}_0(\Omega;\lN)}\\
 &\lesssim_{d,q,\theta,\Omega}
 \|F,\rho f\|_{L^q(\Omega;\lN)}.
 \end{aligned}
\end{equation}

\emph{The case $g=0$ and $1<q<2$.}
Let $G,H$ be smooth test fields, and consider
the adjoint system
\[
 \left\{\begin{aligned}
 -\Delta w_j+\nabla\pi_j+\overline{\lambda_j}w_j
 &=G_j-\div H_j,\\
 \div w_j&=0,\\
 w_j&=0\quad\text{on }\partial\Omega.
 \end{aligned}\right.
\]
Since $q'>2$, the preceding case gives
\[
 \|\rho\nabla w,\rho^2w\|_{L^{q'}(\Omega;\lN)}
 \lesssim_{d,q,\theta,\Omega}
 \|G,\rho H\|_{L^{q'}(\Omega;\lN)}.
\]
Testing the two systems against each other and using
the zero-divergence conditions, we obtain
\begin{equation}\label{eq:duality}
 \begin{aligned}
 \sum_j\int_\Omega
 \rho_j^2u_j\cdot\overline G_j
 +\rho_j^2\nabla u_j:\overline H_j=\sum_j\int_\Omega
 \rho_j^2F_j\cdot\overline w_j
 -\rho_j^2f_j:\nabla\overline w_j.
 \end{aligned}
\end{equation}
By H\"older's inequality and the adjoint estimate,
\[
 \begin{aligned}
 &\left|\sum_j\int_\Omega
 \rho_j^2u_j\cdot\overline G_j
 +\rho_j^2\nabla u_j:\overline H_j\right|\lesssim_{d,q,\theta,\Omega}
 \|F,\rho f\|_{L^q(\Omega;\lN)}
 \|G,\rho H\|_{L^{q'}(\Omega;\lN)}.
 \end{aligned}
\]
Taking the supremum over $G,H$ gives
\eqref{eq:solenoidal-est}. The pressure estimate
\eqref{eq:pressure-global} follows as before.

\emph{The case $g\neq0$.}
Let $b_j=\cB g_j$, where $\cB$ is the Bogovskii
operator from Lemma \ref{lem:Bogovskii}.
By \eqref{eq:Bogovskii},
\[
 \|\rho\nabla b,\rho^2b\|_{L^q(\Omega;\lN)}
 \lesssim_{d,q,\Omega}
 \|\rho^2g\|_{L^q(\Omega;\lN)}.
\]
Set $v_j=u_j-b_j$. Then $\div v_j=0$ and
\[
 -\Delta v_j+\nabla p_j+\lambda_jv_j
 =F_j-\lambda_jb_j+\div(f_j+\nabla b_j).
\]
Since $|\lambda_j|\leq\rho_j^2$ and $\rho_j\geq1$,
\[
 \|F-\lambda b,\rho(f+\nabla b)\|_{L^q(\Omega;\lN)}
 \lesssim_{d,q,\Omega}
 \|F,\rho f,\rho^2g\|_{L^q(\Omega;\lN)}.
\]
Applying the zero-divergence estimate to $(v,p)$
and adding $b$ gives \eqref{eq:general-est}.
This completes the proof.
\end{proof}

\section*{Acknowledgements}
The author thanks Prof.~Zhongwei Shen for suggesting this problem and for many fruitful discussions.

\bibliography{Reference}
\bibliographystyle{abbrv}
\end{document}